\documentclass[12pt]{article}
\usepackage{latexsym,amsfonts,amssymb,amsmath}
\usepackage{color}
\usepackage{colordvi,multicol}

\usepackage{ulem}

\allowdisplaybreaks

\newtheorem{thm}{Theorem}[section]

\newtheorem{lem}[thm]{Lemma}
\newtheorem{pro}[thm]{Proposition}

\newtheorem{rem}[thm]{Remark}

\date{}

\begin{document}

\title{\bf Some New Results on   Gaussian Product Inequalities  }
 \author{ Qian-Qian Zhou, Han Zhao, Ze-Chun Hu\footnote{Corresponding author}, Renming Song\\ \\
 {\small School of Science, Nanjing University of Posts and Telecommunications, Nanjing  210023, China}\\
 {\small qianqzhou@yeah.net}\\
 {\small College of Mathematics, Sichuan  University,
 Chengdu 610065, China}\\
 {\small 1339875802@qq.com}\\
 {\small College of Mathematics, Sichuan  University,
 Chengdu 610065, China}\\
{\small zchu@scu.edu.cn}\\
{\small Department of Mathematics,
University of Illinois, Urbana, IL 61801, USA}\\
 {\small rsong@illinois.edu}}
\maketitle

\begin{abstract}

The long-standing Gaussian product inequality (GPI) conjecture states that, for any centered $\mathbb{R}^n$-valued Gaussian random vector $(X_1, \dots, X_n)$ and any positive reals $\alpha_1, \dots, \alpha_n$, ${\bf E}[\prod_{j=1}^{n}|X_j|^{\alpha_j}]\ge \prod_{j=1}^{n}{\bf E}[|X_j|^{\alpha_j}]$. In this paper, we
present some related inequalities
 for centered $\mathbb{R}^n$-valued Gaussian random vector $(X_1, \dots, X_n)$ when $\{\alpha_1, \dots, \alpha_n\}$ contains both positive and negative numbers.

\end{abstract}



\section{Introduction and main results}

The  Gaussian product inequality (GPI) conjecture can be traced back to \cite{arias}. An earlier form of this conjecture states
that for any centered $\mathbb{R}^n$-valued Gaussian random vector $(X_1, \dots, X_n)$,
\begin{gather}\label{a1}
{\bf E}\left[\prod_{j=1}^{n}X_j^{2m}\right]\ge \prod_{j=1}^n
{\bf E}\left[X_j^{2m}\right],\ m\in \mathbb{N}.
\end{gather}
 Li and Wei  \cite{LW12} proposed the following general version of the GPI conjecture:
 for any centered $\mathbb{R}^n$-valued Gaussian random vector $(X_1, \dots, X_n)$ and any positive reals $\alpha_1, \dots, \alpha_n$,
\begin{gather}\label{a2}
{\bf E}\left[\prod_{j=1}^{n}|X_j|^{\alpha_j}\right]\ge \prod_{j=1}^n{\bf E}\left[|X_j|^{\alpha_j}\right].
\end{gather}
The GPI is known to imply the real polarization conjecture in functional analysis and is related to the U-conjecture,
see \cite{MNPP16}.

The GPI conjecture is still open. However, some special cases have been proved.  Frenkel \cite{F08} proved, using
%ZH an an algebraic method,
an algebraic method,
\eqref{a2}
for $\alpha_j=2$, $j=1, \dots, n$, i.e.,  for any centered  $\mathbb{R}^n$-valued Gaussian random vector $(X_1, \dots, X_n)$,
\begin{gather}\label{a2-1}
{\bf E}\left[\prod_{j=1}^{n}X_j^{2}\right]\ge \prod_{j=1}^n{\bf E}\left[X_j^{2}\right].
\end{gather}
Malicet et al. \cite{MNPP16}, among many other things, gave an analytic proof to (\ref{a2-1}).

It follows from \cite[Corollary 1.1, Theorem 3.1 and Remark 1.4]{KR81} that the GPI \eqref{a2} holds for $n=2$, i.e.,
for any  centered $\mathbb{R}^2$-valued Gaussian random vector $(X_1, X_2)$ and any $\alpha_1, \alpha_2>0$,
 \begin{equation}\label{n=2}
 {\bf E}[|X_1|^{\alpha_1}|X_2|^{\alpha_2}]\ge {\bf E} |X_1|^{\alpha_1}]{\bf E} |X_2|^{\alpha_1}].
 \end{equation}

For the case $n=3$, Lan et al. \cite{LHW20}  proved the following result using
hypergeometric functions:
for any $m_1, m_2\in \mathbb{N}$ and any centered $\mathbb{R}^3$-valued Gaussian random vector $(X_1, X_2, X_3)$,
\begin{eqnarray}\label{1.4}
{\bf E}\left[X_1^{2m_1}X_2^{2m_2}X_3^{2m_2}\right]\ge {\bf E}\left[X_1^{2m_1}\right]{\bf E}\left[X_2^{2m_2}\right]{\bf E}\left[X_3^{2m_2}\right].
\end{eqnarray}
Russell and Sun \cite{RS22-d} proved that for any $m_2, m_3\in \mathbb{N}$ and any centered $\mathbb{R}^3$-valued Gaussian random vector $(X_1, X_2, X_3)$,
\begin{eqnarray}\label{1.5}
{\bf E}\left[X_1^{2}X_2^{2m_2}X_3^{2m_3}\right]\ge {\bf E}\left[X_1^{2}\right]{\bf E}\left[X_2^{2m_2}\right]{\bf E}\left[X_3^{2m_3}\right].
\end{eqnarray}
Herry et al. \cite{HMP22} proved that
for any $m_1,m_2, m_3\in \mathbb{N}$ and any centered $\mathbb{R}^3$-valued Gaussian random vector $(X_1, X_2, X_3)$,
\begin{eqnarray}\label{1.6}
{\bf E}\left[X_1^{2m_1}X_2^{2m_2}X_3^{2m_3}\right]\ge {\bf E}\left[X_1^{2m_1}\right]{\bf E}\left[X_2^{2m_2}\right]{\bf E}\left[X_3^{2m_3}\right].
\end{eqnarray}

Genest and Ouimet \cite{GO22} proved that if there exists a matrix $C\in [0, +\infty)^{n\times n}$ such that\linebreak  $(X_1, X_2, \ldots, X_n)=(Z_1, Z_2, \ldots, Z_n)C$ in law, where $(Z_1, Z_2, \ldots, Z_n)$ is
a standard $\mathbb{R}^n$-valued Gaussian random vector, then
\begin{eqnarray}\label{1.7}
{\bf E}\left[\prod_{j=1}^{n}X_j^{2m_j}\right]\ge \mathbf{E}\left[\prod_{j=1}^kX_j^{2m_j}\right]\mathbf{E}
\left[\prod_{j=k+1}^nX_j^{2m_j}\right],\ m_j\in \mathbb{N},j=1,\ldots,n,\ \forall 1\leq k\leq n-1.
\end{eqnarray}
Russell and Sun \cite{RS22-a} proved, among other things, that \eqref{1.7} holds if all the correlation coefficients of $(X_1, \dots, X_n)$ are nonnegative.
Edelmann et al. \cite{ERR22} extended \eqref{1.7} to the
multi-variate gamma distributions. For other related work,
we refer to Genest and Ouiment \cite {GO22-b},  and  Russell and Sun \cite{RS22-b}.

Wei \cite{W14} proved the following version of GPI: If $(X_1, \dots, X_n)$ is a  centered $\mathbb{R}^n$-valued Gaussian random vector, then  for any $\alpha_1, \dots, \alpha_n\in (-1,0)$,
\begin{eqnarray}\label{a3}
{\bf E}\left[\prod_{j=1}^{n}|X_j|^{\alpha_j}\right]\ge {\bf E}\left[\prod_{j=1}^k|X_j|^{\alpha_j}\right]{\bf E}\left[\prod_{j=k+1}^n|X_j|^{\alpha_j}\right],\ \forall 1\leq k\leq n-1.
\end{eqnarray}
 Russell and Sun \cite{RS22-c}  proved the following opposite  GPI:
 If $(X_1, X_2)$ is a centered $\mathbb{R}^2$-valued Gaussian random vector, then  for any  $\alpha_1\in (-1, 0)$ and $\alpha_2>0$,
 \begin{gather}\label{3}
 {\bf E}[|X_1|^{\alpha_1}|X_2|^{\alpha_2}]\le {\bf E}[|X_1|^{\alpha_1}]{\bf E}[|X_2|^{\alpha_2}].
 \end{gather}
In \cite{HZZ23}, several quantitative versions of the Gaussian product inequality for $n=2$ were given.

In this paper,
we will prove several inequalities related to the GPI \eqref{a2} when $n\ge 3$ and $\{\alpha_1, \dots, \alpha_n\}$
contains both positive and negative numbers.  Our main results are as follows.
	
\begin{thm}\label{th-2}
	Let $\textbf{X}=(X_1,\dots,X_n)$ be
a centered $\mathbb{R}^n$-valued Gaussian random vector with
correlation coefficients $\{\rho_{ij}|1\leq i<j\leq n\}$.
Then for any $\alpha_1\in (0,1)$,
	\begin{align}\label{th-2-a}
		{\bf E}\left[{ \vert X_1 \vert }^{-\alpha_1}\prod_{i=2}^{n}X_i^2\right]
		\geq\prod_{i=2}^{n}(1-\rho_{1i}^{2}){\bf E}[{ \vert X_1 \vert }^{-\alpha_1}]\prod_{i=2}^{n}{\bf E}[X_i^2].
	\end{align}
\end{thm}

\begin{rem}
 (i) By letting $\rho_{1i}=0$  in (\ref{th-2-a}) for any $i=2,\ldots,n$, we know that (\ref{th-2-a}) implies (\ref{a2-1}).

 (ii) For $n=2$, by (\ref{3}) and (\ref{th-2-a}), we know that for any real numbers   $-1<\alpha_1<0,$ $\alpha_2>0$,
 $$
 (1-\rho_{12}^2){\bf E}[|X_1|^{\alpha_1}]{\bf E}[|X_2|^{\alpha_2}]\le {\bf E}[|X_1|^{\alpha_1}|X_2|^{\alpha_2}]\le {\bf E}[|X_1|^{\alpha_1}]{\bf E}[|X_2|^{\alpha_2}]. $$
\end{rem}

\begin{thm}\label{th-1}
	Let $\textbf{X}=(X_1,\dots,X_n)$ be
	a centered $\mathbb{R}^n$-valued Gaussian random vector.
	Then for any  $\alpha_1,\dots, \alpha_{n-1}\in (0,1)$ and $\alpha_n\in (0,\infty)$,
	\begin{align}\label{th-1-a}
{\bf E}\left[{ \vert X_1 \vert }^{-\alpha_1}{\cdots}{ \vert X_{n-1} \vert }^{-\alpha_{n-1}}{ \vert X_n \vert }^{\alpha_n}\right]\leq
{\bf E}\left[\prod_{k=1}^{n-1}{ \vert X_k \vert }^{-\alpha_k}\right]{\bf E}[{ \vert X_n \vert }^{\alpha_n}].
	\end{align}
\end{thm}

\begin{rem}
(i)
By taking $n=2$  in \eqref{th-1-a}, we get \eqref{3}.

(ii)
If the covariance matrix $\Sigma$ of $(X_1, \dots, X_n)$ can be written as $MM'$, where $M$ is a lower triangular matrix and $M'$ stands for the transpose of $M$, then combining (\ref{th-1-a}) and  Wei \cite[Theorem 3.1]{W14}, we get that
for any  $\alpha_1,\dots, \alpha_{n-1}\in (0,1)$ and $\alpha_n\in (0,\infty)$,
\begin{gather}\label{2}
{\bf E}\left[{ \vert X_1 \vert }^{-\alpha_1}{\cdots}{ \vert X_{n-1} \vert }^{-\alpha_{n-1}}{ \vert X_n \vert }^{\alpha_n}\right]\leq
\prod_{j=1}^{n-1} \left(\frac{1}{m_{j,j}}\right)^{\alpha_j}\, \prod_{j=1}^{n-1} {\bf E}[|X_j|^{-\alpha_j}] {\bf E}[|X_n|].
\end{gather}
\end{rem}

\begin{pro}\label{co-1}
	Let $\textbf{X}=(X_1,\dots,X_4)$ be  a centered  $\mathbb{R}^4$-valued Gaussian random vector with correlation coefficients $\{\rho_{ij}|1\leq i<j\leq 4\}$.
Then for any $\alpha_1\in (0,1)$ and
	any positive even integers $ \alpha_2, \alpha_3, \alpha_4$,
	\begin{align}\label{co-1-a}
		{\bf E}\left[{ \vert X_1 \vert }^{-\alpha_1}{ \vert X_2 \vert }^{\alpha_2}{ \vert X_3 \vert }^{\alpha_3}{ \vert X_4 \vert }^{\alpha_4}\right]
		\geq\prod_{i=2}^{4}(1-\rho_{1i}^{2})^{\frac{\alpha_i}{2}}{\bf E}[{ \vert X_1 \vert }^{-\alpha_1}]\prod_{i=2}^{4}{\bf E}[{ \vert X_i \vert }^{\alpha_i}].
	\end{align}
\end{pro}

By letting $\rho_{1i}=0$  in \eqref{co-1-a} for any
$i=2,3,4,$  we easily get \eqref{1.6}.

\begin{pro}\label{th-3}
	Let $\textbf{X}=(X_1,X_2,X_3)$ be  a centered  $\mathbb{R}^3$-valued Gaussian random vector with correlation coefficients $\{\rho_{ij}|1\leq i<j\leq 3\}$.
Then for any  $\alpha_1\in (0,1)$ and $\alpha_2,\alpha_3\in (0,\infty)$,
	\begin{eqnarray}\label{th-3-a}
\prod_{i=2}^{3} (1-\rho_{1i}^{2})^{\frac{\alpha_i}{2}}{\bf E}[{ \vert X_1 \vert }^{-\alpha_1}]{\bf E}[{ \vert X_{2} \vert }^{\alpha_{2}}]{\bf E}[{ \vert X_3 \vert }^{\alpha_3}]
		&\leq&{\bf E}\left[{ \vert X_1 \vert }^{-\alpha_1}{ \vert X_{2} \vert }^{\alpha_{2}}{ \vert X_3 \vert }^{\alpha_3}\right]\nonumber\\
		&\leq&\textbf{C}_{\alpha_2,\alpha_3}{\bf E}[{ \vert X_1 \vert }^{-\alpha_1}]{\bf E}[{ \vert X_{2} \vert }^{\alpha_{2}}]{\bf E}[{ \vert X_3 \vert }^{\alpha_3}],\quad\quad
	\end{eqnarray}
	where
$ \textbf{C}_{\alpha_2,\alpha_3}=\frac{\Gamma(\frac{1}{2})\Gamma(\frac{1}{2}+
\frac{\alpha_2+\alpha_3}{2})}{\Gamma(\frac{1}{2}+\frac{\alpha_2}{2})
\Gamma(\frac{1}{2}+\frac{\alpha_3}{2})} $.
\end{pro}

\begin{pro}\label{th-5}
	Let $\textbf{X}=(X_1,X_2,X_3)$ be a  centered  $\mathbb{R}^3$-valued Gaussian random vector with covariance matrix $\Sigma=\begin{pmatrix}
		1& a& b\\
		a& 1& c\\
		b& c& 1
	\end{pmatrix}$. Then for any $\alpha_1,\alpha_2\in (0,1)$ and $\alpha_3\in (0,\infty)$,
	\begin{align}\label{th-5-a}
		\frac{\det(\Sigma)}{1-a^{2}}{\bf E}[{ \vert X_1 \vert }^{-\alpha_1}]{\bf E}[{ \vert X_{2} \vert }^{-\alpha_{2}}]{\bf E}[{ \vert X_3 \vert }^{\alpha_3}]
		\leq{\bf E}[{ \vert X_1 \vert }^{-\alpha_1}{ \vert X_{2} \vert }^{-\alpha_{2}}{ \vert X_3 \vert }^{\alpha_3}]
		\leq{\bf E}[{ \vert X_1 \vert }^{-\alpha_1}{ \vert X_{2} \vert }^{-\alpha_{2}}]{\bf E}[{ \vert X_3 \vert }^{\alpha_3}].
	\end{align}
\end{pro}

\section{Proofs }\setcounter{equation}{0}

In this section, we give the
%ZH poofs
proofs
of our main results. In Subsection 2.1, we first recall some lemmas that will be used in proving our main results.

\subsection{Preliminaries}

\begin{lem}
(Bernstein \cite[Fact 2.8.16]{BD09})\label{lem_1}
Let 	$ A $  be an $ m\times m $ invertible matrix, $ B \in R^{m\times n}$, $ C \in R^{n\times m}$, $ D \in R^{n\times n}$, and $ (D-CA^{-1}B) $ be an $ n\times n $ invertible matrix. Then
	\begin{align*}
		\begin{pmatrix}
			A &B\\
			C &D
		\end{pmatrix}^{-1}=\begin{pmatrix}
			A^{-1}+A^{-1}B(D-CA^{-1}B)^{-1}CA^{-1}&-A^{-1}B(D-CA^{-1}B)^{-1}\\
			-(D-CA^{-1}B)^{-1}CA^{-1}&(D-CA^{-1}B)^{-1}
		\end{pmatrix}.
	\end{align*}
\end{lem}

\begin{lem}
(Bernstein \cite[Fact 2.8.17]{BD09})\label{lem_2}
Let $ A $  be an $ m\times m $ invertible matrix, $ B \in R^{m\times n}$, $ C \in R^{n\times m}$, $ D \in R^{n\times n}$, and $ (A-BD^{-1}C) $ be an $ m\times m $ invertible matrix. Then
	\begin{align*}
		\begin{pmatrix}
			A &B\\
			C &D
		\end{pmatrix}^{-1}=\begin{pmatrix}
			(A-BD^{-1}C)^{-1}&-(A-BD^{-1}C)^{-1}BD^{-1}\\
			-D^{-1}C(A-BD^{-1}C)^{-1}&D^{-1}+D^{-1}C(A-BD^{-1}C)^{-1}BD^{-1}
		\end{pmatrix}.
	\end{align*}
\end{lem}

\begin{lem}\label{lem-3}
Suppose that $\Sigma=(\sigma_{ij})_{n \times n}=\begin{pmatrix}
		\Sigma_{1} & \textbf{t}\\
		\textbf{t}^{'}          &1
	\end{pmatrix}$
is a positive definite
$n\times n$ matrix with all diagonal elements equal to 1, where
$ \Sigma_{1}=(\sigma_{ij})_{(n-1) \times (n-1)} $
is an $(n-1)\times (n-1)$ matrix and
$\textbf{t}=(\sigma_{1n},\dots,\sigma_{n-1\  n})^{'}$,
then  $ \Sigma_{1}-\textbf{t}\textbf{t}^{'} $ is also a positive definite matrix.
\end{lem}
{\bf Proof.}
We prove this lemma by induction.
When $n=2$, we have
$ \Sigma=\begin{pmatrix}
	\Sigma_{1} & \sigma_{12}\\
	\sigma_{21}          &1
\end{pmatrix} $,
where $ \Sigma_{1}=1 $ and
$\sigma_{12}=\sigma_{21} $. Thus $ \Sigma_{1}-tt^{'}=1-\sigma_{12}^{2}\textgreater0 $
by the positive definiteness of $\Sigma$.
Thus, the result holds for $n=2.$

Suppose that the result holds for $n-1 (n\geq 3)$: If  $\Sigma=\begin{pmatrix}
	\Sigma_{1} & \textbf{t}_{n-2}\\
	\textbf{t}_{n-2}^{'}          &1
\end{pmatrix}_{(n-1)\times(n-1)}$
is a  positive definite $(n-1)\times (n-1)$ matrix with all diagonal elements equal to 1,
where $ \Sigma_{1}=(\sigma_{ij})_{(n-2)\times(n-2)}$
is an $(n-2)\times (n-2)$ matrix and
$\textbf{t}_{n-2}=(\sigma_{1\,n-1},\dots,\sigma_{n-2\,  n-1})^{'} $,
then ${\Sigma_{1}}-\textbf{t}_{n-2}\textbf{t}_{n-2}^{'}$ is also a positive definite matrix.
We will prove that the result holds for $n$. Assuming that $\Sigma=\begin{pmatrix}
	\Sigma_{1}& \textbf{t}_{n-1}\\
	\textbf{t}_{n-1}^{'}          &1
\end{pmatrix}_{n\times n}$  is a
positive definite  $n\times n$ matrix with all diagonal elements equal to 1, where
$\Sigma_{1}=(\sigma_{ij})_{(n-1)\times(n-1)} $
is an $(n-1)\times (n-1)$ matrix and
$\textbf{t}_{n-1}=(\sigma_{1n},\dots,\sigma_{n-1\ n})^{'} $.
It follows from \cite[Fact 2.14.3]{BD09} that $\det({\Sigma_{1}}-\textbf{t}_{n-1}\textbf{t}_{n-1}^{'})=\det(\Sigma)\textgreater0$.
Thus
we only need to prove that the $(n-2)$-th order principal submatrix  of $ ({\Sigma_{1}}-\textbf{t}_{n-1}\textbf{t}_{n-1}^{'}) $  is positively definite.

Note that
\begin{align*}
	{\Sigma_{1}}-\textbf{t}_{n-1}\textbf{t}_{n-1}^{'}=\begin{pmatrix}
		\widetilde{\Sigma}_1
		& \widetilde{{\bf t}}_{n-2}\\
		\widetilde{{\bf t}}^{'}_{n-2} & 1
	\end{pmatrix}-\begin{pmatrix}
		\hat{\textbf{t}}_{n-2}\\
	\sigma_{n-1\  n}
	\end{pmatrix}\begin{pmatrix}
        \hat{\textbf{t}}_{n-2}^{'} & \sigma_{n\  n-1}
	\end{pmatrix},
\end{align*}
where $\widetilde{\Sigma}_{1}={(\sigma_{ij})}_{(n-2)\times(n-2)}$,
$\widetilde{{\bf t}}_{n-2}=(\sigma_{1\ n-1}, \ldots, \sigma_{n-2\ n-1})^{'}$, $\hat{\textbf{t}}_{n-2}=(\sigma_{1n},\dots,\sigma_{n-2\ n})^{'}$.
So the $(n-2)$-th order principal submatrix  of $ ({\Sigma_{1}}-\textbf{t}_{n-1}\textbf{t}_{n-1}^{'}) $  is
$({\widetilde{\Sigma}_{1}}- \hat{\textbf{t}}_{n-2}\hat{\textbf{t}}_{n-2}^{'})$.
Note that $ \begin{pmatrix}
	{\widetilde{\Sigma}_{1}} &\hat{\textbf{t}}_{n-2}\\
	\hat{\textbf{t}}_{n-2}^{'} &1
\end{pmatrix} $ is the $(n-1)$-th order principal submatrix  of $\Sigma$.
Since $\Sigma$ is positive definite, $ \begin{pmatrix}
	{\widetilde{\Sigma}_{1}} &\hat{\textbf{t}}_{n-2}\\
	\hat{\textbf{t}}_{n-2}^{'} &1
\end{pmatrix} $ is also a positive definite   matrix.
Now using the induction hypothesis, we get that $ ({\widetilde{\Sigma}_{1}}- \hat{\textbf{t}}_{n-2}\hat{\textbf{t}}_{n-2}^{'})$ is a positive definite matrix. The proof is complete.
\hfill\fbox

\subsection{Proof of Theorem   \ref{th-2}}

 Without loss of generality, we can assume that
 $X_1, \dots, X_n$ are all standard Gaussian random variables.
Then, the desired assertion reduces to
\begin{align}\label{th-2-a-rs}
		{\bf E}\left[{ \vert X_1 \vert }^{-\alpha_1}\prod_{i=2}^{n}X_i^2\right]
		\geq\prod_{i=2}^{n}(1-\rho_{1i}^{2}){\bf E}[{ \vert X_1 \vert }^{-\alpha_1}].
	\end{align}

Without loss of generality, we can assume that the covariance matrix $\Sigma$ is
invertible and thus positive definite.
Otherwise, we can take the covariance matrix to be $\frac{1}{1+1/m}(\Sigma+\frac{1}{m}I_n)$ and then let $m\to\infty$. Hereafter, $I_n$ stands for the $n$-dimensional identity matrix.

For any $t>0$, let  $T$ be the $n\times n$ diagonal matrix $\mbox{diag}(t^{\frac{2}{\alpha_1}},0,\dots,0)$. Since $\Sigma$ is positive definite, we know that $\Sigma^{-1}$ is also positive definite and so
$\Sigma^{-1}+2T$ is also positive definite. Let $(Y_1, \dots, Y_n)$ be a centered $\mathbb{R}^n$-valued Gaussian random vector with covariance matrix $(\Sigma^{-1}+2T)^{-1}$ and let
	\begin{align*}
		\hat{f}(x_2,\dots,x_n):=\int_{\mathbb{R}}\frac{\exp(-\frac{1}{2}\left\langle \textbf{x},(\Sigma^{-1}+2T)\textbf{x}\right\rangle)}{(2\pi)^{\frac{n}{2}}
\sqrt{\det((\Sigma^{-1}+2T)^{-1})}}dx_1,
	\end{align*}
where $\textbf{x}=(x_1, x_2, \ldots, x_n)$.

By the definition of the Gamma function and a change of variables, we have
\begin{eqnarray}\label{p-th-2-a}
		{ \vert u \vert }^{-\alpha}=\frac{1}{\Gamma(1+\frac{\alpha}{2})}\int_{0}^{\infty}
\exp\left(-t^{\frac{2}{\alpha}}u^{2}\right)dt,\
		\forall \alpha>0.
\end{eqnarray}
By (\ref{p-th-2-a}),  Fubini's theorem and (\ref{a2-1}), we  have
\begin{align}\label{p-th-2-b}
&{\bf E}\left[{ \vert X_1 \vert }^{-\alpha_1}\prod_{i=2}^{n}
X_i^2\right]\nonumber\\
&= \frac{1}{\Gamma(1+\frac{\alpha_{1}}{2})} \int_0^{\infty} {\bf E}\left[\exp\left(-t^{\frac{2}{\alpha_1}}X^2_1\right) \prod^{n}_{i=2}X_i^2  \right]dt      \nonumber\\
&= \frac{1}{\Gamma(1+\frac{\alpha_{1}}{2})} \int_0^{\infty} \int_{\mathbb{R}^n} \prod_{i=2}^n x_i^2 \exp(-\left\langle \textbf{x},T\textbf{x} \right\rangle)\frac{\exp(-\frac{1}{2}\left\langle \textbf{x},\Sigma^{-1} \textbf{x} \right\rangle)}{(2\pi)^{\frac{n}{2}}\sqrt{\det(\Sigma)}}d\textbf{x}dt   \nonumber\\
&=  \frac{1}{\Gamma(1+\frac{\alpha_{1}}{2})} \int_0^{\infty} \left(\int_{\mathbb{R}^n} \prod_{i=2}^n x_i^2 \frac{\exp(-\frac{1}{2}\left\langle \textbf{x},(\Sigma^{-1}+2T)\textbf{x}\right\rangle)}{(2\pi)^{\frac{n}{2}}
\sqrt{\det((\Sigma^{-1}+2T)^{-1})}}d\textbf{x}\right)
\frac{\sqrt{\det((\Sigma^{-1}+2T)^{-1})}}
{\sqrt{\det(\Sigma)}}dt\nonumber\\	
	&=\frac{1}{\Gamma(1+\frac{\alpha_{1}}{2})}
\int_{0}^{\infty}[\det(I_n+2T\Sigma)]^{-\frac{1}{2}}\int_{\mathbb{R}^{n}}\prod_{i=2}^{n}
x_i^2\frac{\exp(-\frac{1}{2}\left\langle \textbf{x},(\Sigma^{-1}+2T)\textbf{x}\right\rangle)}{(2\pi)^{\frac{n}{2}}
\sqrt{\det((\Sigma^{-1}+2T)^{-1})}}d\textbf{x}dt\nonumber\\
		&=\frac{1}{\Gamma(1+\frac{\alpha_{1}}{2})}
\int_{0}^{\infty}[\det(I_n+2T\Sigma)]^{-\frac{1}{2}}\int_{\mathbb{R}^{n-1}}\prod_{i=2}^{n}
x_i^2\hat{f}(x_2,\dots,x_n)dx_{2}\dots dx_{n}dt\nonumber\\
&=\frac{1}{\Gamma(1+\frac{\alpha_{1}}{2})}\int_{0}^{\infty}
\left(1+2t^{\frac{2}{\alpha_1}}\right)^{-\frac{1}{2}}
{\bf E}\left[\prod_{i=2}^{n}Y_i^2\right]dt\nonumber\\
		&\geq\frac{1}{\Gamma(1+\frac{\alpha_{1}}{2})}\int_{0}^{\infty}
\left(1+2t^{\frac{2}{\alpha_1}}\right)^{-\frac{1}{2}}\prod_{i=2}^{n}	
	{\bf E}[Y_i^2]dt\nonumber\\
&=\frac{1}{\Gamma(1+\frac{\alpha_{1}}{2})}\int_{0}^{\infty}
\left(1+2t^{\frac{2}{\alpha_1}}\right)^{-\frac{1}{2}}\prod_{i=2}^{n}
{\rm Var}(Y_i)dt.
	\end{align}

The covariance matrix $ (\Sigma^{-1}+2T)^{-1}$ of $(Y_1,\dots,Y_n)$ can be written as $\begin{pmatrix}
a& \mathbf{b}\\
\mathbf{b}' & D
	\end{pmatrix}$, where $D$
	is an $(n-1)\times (n-1)$ matrix and $ \prod_{i=2}^{n}{\rm Var}(Y_i) $ is the product of  the diagonal elements  of the  matrix $D.$
We rewrite the matrix $\Sigma$ as $\begin{pmatrix}
		1 & {\bf t}\\
		{\bf t}^{'} &\Sigma_{1}
	\end{pmatrix}$
 with $ \Sigma_{1} $ being an $(n-1)\times (n-1)$ matrix and $ {\bf t}=(\rho_{12},\dots,\rho_{1n})$.
 By the property of determinants of block matrices,
 we have $\det(\Sigma)=\det(\Sigma_{1}-{\bf t}^{'}{\bf t})>0$, which implies that  $\Sigma_{1}-{\bf t}^{'}{\bf t}$ is invertible. Then,  by Lemma \ref{lem_1}, we have
	\begin{eqnarray}\label{p-th-2-c}
		\Sigma^{-1}+2T=\begin{pmatrix}
			1+{\bf t}(\Sigma_{1}-{\bf t}^{'}{\bf t})^{-1}{\bf t}^{'}+2t^{\frac{2}{\alpha_1}} & -{\bf t}(\Sigma_{1}-{\bf t}^{'}{\bf t})^{-1}\\
			-(\Sigma_{1}-{\bf t}^{'}{\bf t})^{-1}{\bf t}^{'} & (\Sigma_{1}-{\bf t}^{'}{\bf t})^{-1}
		\end{pmatrix}.
	\end{eqnarray}
Since $\Sigma^{-1}+2T$ is positive definite, by (\ref{p-th-2-c}) and Lemma \ref{lem_2}, we have
	\begin{align*}
		D&=(\Sigma_{1}-{\bf t}^{'}{\bf t})+{\bf t}^{'}\left(1+{\bf t}(\Sigma_{1}-{\bf t}^{'}{\bf t})^{-1}{\bf t}^{'}+2t^{\frac{2}{\alpha_1}}-{\bf t}(\Sigma_{1}-{\bf t}^{'}{\bf t})^{-1}(\Sigma_{1}-{\bf t}^{'}{\bf t})(\Sigma_{1}-{\bf t}^{'}{\bf t})^{-1}{\bf t}^{'}\right)^{-1}{\bf t}\\
		&=(\Sigma_{1}-{\bf t}^{'}{\bf t})+\left(1+2t^{\frac{2}{\alpha_1}}\right)^{-1}{\bf t}^{'}{\bf t}\\
		&=\Sigma_{1}-\left(1-\frac{1}{1+2t^{\frac{2}{\alpha_1}}}\right){\bf t}^{'}{\bf t}.
	\end{align*}
It follows that
\begin{eqnarray}\label{p-th-2-d}
 \prod_{i=2}^{n}{\rm Var}(Y_i)
 =\prod_{i=2}^{n}\left[1-\left(1-\frac{1}{1+2t^{\frac{2}{\alpha_1}}}\right)
\rho_{1i}^{2}\right]\geq\prod_{i=2}^{n}(1-\rho_{1i}^{2}).
\end{eqnarray}
Then, by (\ref{p-th-2-b}), (\ref{p-th-2-d}) and (\ref{p-th-2-a}), we obtain
	\begin{align*}
		{\bf E}\left[{ \vert X_1 \vert}^{-\alpha_1}\prod_{i=2}^{n}{ \vert X_i \vert }^{2}\right]
		&=\frac1{\Gamma(1+\frac{\alpha_{1}}{2})}
\int_{0}^{\infty}\left(1+2t^{\frac{2}{\alpha_1}}\right)^{-\frac{1}{2}}
\prod_{i=2}^{n}{\rm Var}(Y_i)dt\\
		&\geq
		\frac{\prod_{i=2}^{n}(1-\rho_{1i}^{2})}{\Gamma(1+\frac{\alpha_{1}}{2})}\int_{0}^{\infty}
(1+2t^{\frac{2}{\alpha_1}})^{-\frac{1}{2}}dt\\
		&=\prod_{i=2}^{n}(1-\rho_{1i}^{2})
		{\bf E}[{ \vert X_1 \vert }^{-\alpha_1}],
	\end{align*}
which implies that \eqref{th-2-a-rs} holds. The proof is complete.

\subsection{Proof of Theorem   \ref{th-1} }

Without loss of generality, we can assume $X_1, \dots, X_n$ are all standard Gaussian random variables and that the covariance matrix $\Sigma$ is  invertible and thus positive definite. For any $t_1, t_2, \dots, t_{n-1}>0$, let  $T$ be the $n\times n$ diagonal matrix $\mbox{diag}(t_1^{\frac{2}{\alpha_1}},\dots, t_{n-1}^{\frac{2}{\alpha_{n-1}}}, 0)$. Since $\Sigma$ is positive definite, we know that $\Sigma^{-1}$ is also positive definite and so
$\Sigma^{-1}+2T$ is also positive definite. Let $(Y_1, \dots, Y_n)$ be a centered $\mathbb{R}^n$-valued Gaussian random vector with covariance matrix $(\Sigma^{-1}+2T)^{-1}$ and let
	\begin{align*}
		\hat{f}(x_n):=\int_{\mathbb{R}^{n-1}}\frac{\exp(-\frac{1}{2}\left\langle \textbf{x},(\Sigma^{-1}+2T)\textbf{x}\right\rangle)}{(2\pi)^{\frac{n}{2}}\sqrt{\det((\Sigma^{-1}+2T)^{-1})}}dx_1\dots dx_{n-1},
	\end{align*}
 where $\textbf{x}=(x_1, x_2, \ldots, x_n)$.

Similar to the proof of Theorem  \ref{th-2}, we have that
\begin{eqnarray}\label{proo-th-1-a}
	\hspace{-0.5cm}	&&{\bf E}\left[{ \vert X_1 \vert }^{-\alpha_1}{\dots}{ \vert X_{n-1} \vert }^{-\alpha_{n-1}}{ \vert X_n \vert }^{\alpha_n}\right]\nonumber\\		\hspace{-0.5cm}&&=\frac{1}{\prod_{i=1}^{n-1}\Gamma(1+\frac{\alpha_i}{2})}\int_{0}^{\infty}\cdots
\int_{0}^{\infty}\int_{\mathbb{R}^{n}}[\det(I_n+2T\Sigma)]^{-\frac{1}{2}}{ \vert x_n \vert }^{\alpha_n}\frac{\exp(-\frac{1}{2}\left\langle \textbf{x},(\Sigma^{-1}+2T)\textbf{x}\right\rangle)}{(2\pi)^{\frac{n}{2}}
\sqrt{\det((\Sigma^{-1}+2T)^{-1})}}d\textbf{x}dt_{1}\dots dt_{n-1}\nonumber\\		\hspace{-0.5cm}&&=\frac{1}{\prod_{i=1}^{n-1}\Gamma(1+\frac{\alpha_i}{2})}\int_{0}^{\infty}\cdots
\int_{0}^{\infty}[\det(I_n+2T\Sigma)]^{-\frac{1}{2}}\int_{\mathbb{R}}{ \vert x_n \vert }^{\alpha_n}\hat{f}(x_{n})dx_{n}dt_{1}\dots dt_{n-1}\nonumber\\
\hspace{-0.5cm}		&&=\frac{1}{\prod_{i=1}^{n-1}\Gamma(1+\frac{\alpha_i}{2})}\int_{0}^{\infty}\cdots
\int_{0}^{\infty}[\det(I_n+2T\Sigma)]^{-\frac{1}{2}}
{\bf E}[{ \vert Y_n \vert }^{\alpha_n}]dt_{1}\dots dt_{n-1}\nonumber\\
\hspace{-0.5cm}		&&=\frac{1}{\prod_{i=1}^{n-1}\Gamma(1+\frac{\alpha_i}{2})}\int_{0}^{\infty}\cdots
\int_{0}^{\infty}[\det(I_n+2T\Sigma)]^{-\frac{1}{2}}{\bf E}[{ \vert X_n \vert }^{\alpha_n}]
({\rm Var}(Y_n))^{\frac{\alpha_n}{2}}dt_{1}\dots dt_{n-1},
\end{eqnarray}
where we used the fact that
\begin{eqnarray}\label{proo-th-1-b}
{\bf E}[{ \vert Y_n \vert }^{\alpha_n}]
&=&\int_{-\infty}^{\infty}|u|^{\alpha_n}\frac{1}
{ \sqrt{2\pi{\rm Var}(Y_n)}}
e^{-\frac{u^2}{2{\rm Var}(Y_n)}}du\nonumber\\
&=&\frac{2^{\frac{\alpha_n}{2}}}{\sqrt{\pi}}\Gamma\left(\frac{\alpha_n+1}{2}\right)
({\rm Var}(Y_n))^{\frac{\alpha_n}{2}}\nonumber\\
&=&{\bf E}[{ \vert X_n \vert }^{\alpha_n}]
({\rm Var}(Y_n))^{\frac{\alpha_n}{2}}.
\end{eqnarray}
We know that ${\rm Var}(Y_n)$ is the last  diagonal element   of
the   matrix
 $ (\Sigma^{-1}+2T)^{-1}$.

We write $\Sigma$ and $T$ as $\Sigma=\begin{pmatrix}
		\Sigma_{1} & {\bf t}\\
		{\bf t}^{'}          &1
	\end{pmatrix}$ and $ T=\begin{pmatrix}
		T_{1} & {\bf 0}\\
		{\bf 0}     & 0
	\end{pmatrix} $ with $ \Sigma_{1} $ and $ T_{1} $
	$(n-1)\times (n-1)$ matrices and
	$ {\bf t}=(\rho_{1 n},\dots,\rho_{n-1\ n})^{'}$.  By Lemma \ref{lem-3}, we know $ \Sigma_{1}-{\bf tt}^{'} $ is a positive definite matrix and thus  invertible. Then,  by  Lemma \ref{lem_2},  we get  that
	\begin{align*}
		\Sigma^{-1}+2T=\begin{pmatrix}
			(\Sigma_{1}-{\bf tt}^{'})^{-1}+2T_{1} & -(\Sigma_{1}-{\bf tt}^{'})^{-1}{\bf t}\\
			-{\bf t}^{'}(\Sigma_{1}-{\bf tt}^{'})^{-1} & 1+{\bf t}^{'}(\Sigma_{1}-{\bf tt}^{'})^{-1}{\bf t}
		\end{pmatrix}.
	\end{align*}
Since the matrix $\Sigma^{-1}+2T$  is positive definite, by  Lemma \ref{lem_1},
we know  that
\begin{eqnarray}\label{proo-th-1-c}
		 {\rm Var}(Y_n)
		&=&\left[1+{\bf t}^{'}(\Sigma_{1}-{\bf tt}^{'})^{-1}{\bf t}-{\bf t}^{'}(\Sigma_{1}-{\bf tt}^{'})^{-1}[(\Sigma_{1}-{\bf tt}^{'})^{-1}+2T_{1}]^{-1}(\Sigma_{1}-{\bf tt}^{'})^{-1}{\bf t}\right]^{-1}\nonumber\\
		&=&\left[1+{\bf t}^{'}(\Sigma_{1}-{\bf tt}^{'})^{-1}{\bf t}-{\bf t}^{'}(\Sigma_{1}-{\bf tt}^{'})^{-1}[I_{n-1}+2(\Sigma_{1}-{\bf tt}^{'})T_{1}]^{-1}{\bf t}\right]^{-1}\nonumber\\
		&=&\left[1+{\bf t}^{'}(\Sigma_{1}-{\bf tt}^{'})^{-1}\left[I_{n-1}-[I_{n-1}+2(\Sigma_{1}-{\bf tt}^{'})T_{1}]^{-1}\right]{\bf t}\right]^{-1}\nonumber\\
		&=&\left[1+{\bf t}^{'}(\Sigma_{1}-{\bf tt}^{'})^{-1}[I_{n-1}+2(\Sigma_{1}-{\bf tt}^{'})T_{1}]^{-1}[2(\Sigma_{1}-{\bf tt}^{'})T_{1}]{\bf t}\right]^{-1}\nonumber\\
		&=&\left[1+2{\bf t}^{'}[(\Sigma_{1}-{\bf tt}^{'})+2(\Sigma_{1}-{\bf tt}^{'})T_{1}(\Sigma_{1}-{\bf tt}^{'})]^{-1}(\Sigma_{1}-{\bf tt}^{'})T_{1}{\bf t}\right]^{-1}\nonumber\\
		&=&\left[1+2{\bf t}^{'}[I_{n-1}+2T_{1}(\Sigma_{1}-{\bf tt}^{'})]^{-1}T_{1}{\bf t}\right]^{-1}\nonumber\\
		&=&\left[1+2{\bf t}^{'}[T_{1}^{-1}+2(\Sigma_{1}-{\bf tt}^{'})]^{-1}{\bf t}\right]^{-1}\le 1,
	\end{eqnarray}
where we used the fact that  $ \Sigma_{1}-{\bf tt}^{'} $ is a positive definite matrix.

According to Sylvester's determinant theorem (see \cite{H85}),
\begin{eqnarray}\label{proo-th-1-d}
		\det(I_n+2T\Sigma)&=&\det(I_n+2T^{\frac{1}{2}}T^{\frac{1}{2}}\Sigma)=
\det(I_n+2T^{\frac{1}{2}}\Sigma T^{\frac{1}{2}})\nonumber\\
		&=&\det\begin{pmatrix}
			I_{n-1}+2T_1^{\frac{1}{2}}\Sigma_{1}T_1^{\frac{1}{2}} &{\bf 0}\\
			{\bf 0}& 1
		\end{pmatrix}\nonumber\\
&=&\det(I_{n-1}+2T_1^{\frac{1}{2}}\Sigma_{1}T_1^{\frac{1}{2}})\nonumber\\
		&=&\det(I_{n-1}+2T_1\Sigma_{1}).
	\end{eqnarray}

By (\ref{proo-th-1-a})-(\ref{proo-th-1-d}),  we have
	\begin{align*}
&		{\bf E}\left[{ \vert X_1 \vert }^{-\alpha_1}{\dots}{ \vert X_{n-1} \vert }^{-\alpha_{n-1}}{ \vert X_n \vert }^{\alpha_n}\right]\\		&\leq\frac{1}{\prod_{i=1}^{n-1}\Gamma(1+\frac{\alpha_i}{2})}\int_{0}^{\infty}\cdots\int_{0}^{\infty}[\det(I_n+2T\Sigma)]^{-\frac{1}{2}}{\bf E}[{ \vert X_n \vert }^{\alpha_n}]dt_{1}\dots dt_{n-1}\\
		&={\bf E}[{ \vert X_n \vert }^{\alpha_n}]\cdot\frac{1}{\prod_{i=1}^{n-1}\Gamma(1+\frac{\alpha_i}{2})}\int_{0}^{\infty}\cdots\int_{0}^{\infty}[\det(I_{n-1}+2T_1\Sigma_{1})]^{-\frac{1}{2}}dt_{1}\dots dt_{n-1}\\
		&={\bf E}\left[{ \vert X_1 \vert }^{-\alpha_1}{\dots}{ \vert X_{n-1}\vert}^{-\alpha_{n-1}}\right]{\bf E}[{\vert X_n \vert }^{\alpha_n}].
	\end{align*}
The proof is complete.

\subsection{Proof of Proposition \ref{co-1}}

 We need only to replace
 $ {\bf E}\left[\prod_{i=2}^{n}{ \vert Y_i \vert }^{2}\right]\geq \prod_{i=2}^{n}{\bf E}\left[{ \vert Y_i \vert }^{2}\right]$ by
$$
 {\bf E}[{ \vert Y_2 \vert }^{\alpha_2}{ \vert Y_3 \vert }^{\alpha_3}{ \vert Y_4 \vert }^{\alpha_4}]\geq {\bf E}[{ \vert Y_2 \vert }^{\alpha_2}]{\bf E}{ \vert Y_3 \vert }^{\alpha_3}]{\bf E}[{ \vert Y_4 \vert }^{\alpha_4}]
 $$
 (which is valid by \eqref{1.6}) in the proof of Theorem \ref{th-2}.
 We omit the details.

\subsection{Proof of Proposition \ref{th-3}}

Letting $n=3$ and replacing $|X_2|^2 |X_3|^2$  by $|X_2|^{\alpha_2}|X_3|^{\alpha_3}$ in the proof of Theorem \ref{th-2}, we
get the  inequality on the left side of \eqref{th-3-a}.
Now we prove the other inequality in \eqref{th-3-a}.

For $t>0$, let $T:={\rm diag}(t^{\frac{2}{\alpha_1}},0,0)$. Then $\Sigma^{-1}+2T$ is positive definite. Let $(Y_1, Y_2, Y_3)$
be a
%ZH centerd
centered
$\mathbb{R}^3$-valued Gaussian random vector with covariance matrix $( \Sigma^{-1}+2T)^{-1}$ and let
\begin{align*}
		\hat{f}(x_2,x_3):=\int_{\mathbb{R}}\frac{\exp(-\frac{1}{2}\left\langle \textbf{x},(\Sigma^{-1}+2T)\textbf{x}\right\rangle)}{(2\pi)^{\frac{3}{2}}\sqrt{\det((\Sigma^{-1}+2T)^{-1})}}dx_1,
	\end{align*}
where  $\textbf{x}=(x_1, x_2, x_3)$.
Similar to the proof of Theorem \ref{th-2}, we  have
	\begin{eqnarray}\label{proof-th-2-a}
		{\bf E}[{ \vert X_1 \vert }^{-\alpha_1}{ \vert X_{2} \vert }^{\alpha_{2}}{ \vert X_3 \vert }^{\alpha_3}]
		=\frac{1}{\Gamma(1+\frac{\alpha_{1}}{2})}\int_{0}^{\infty}
\left(1+2t^{\frac{2}{\alpha_1}}\right)^{-\frac{1}{2}}
 {\bf E}[{ \vert Y_2 \vert }^{\alpha_2}{ \vert Y_3 \vert }^{\alpha_3}]dt,
	\end{eqnarray}
and ${\rm Var}(Y_2){\rm Var}(Y_3)=\prod_{i=2}^{3}
\left[1-\left(1-\frac{1}{1+2t^{\frac{2}{\alpha_1}}}\right)\rho_{1i}^{2}\right]\leq1$.

	By Nabeya \cite{NS51}, we have that
	\begin{eqnarray}\label{proof-th-2-b}
		{\bf E}[{ \vert Y_2 \vert }^{\alpha_2}{ \vert Y_3 \vert }^{\alpha_3}]
		&=&
		{\bf E}[{ \vert Y_2 \vert }^{\alpha_2}]{\bf E}[{ \vert Y_3 \vert }^{\alpha_3}]
		\,F\left(-\frac{\alpha_2}{2},-\frac{\alpha_3}{2};\frac{1}{2};\rho_{t}^{2}\right)\nonumber\\
		&=&{\bf E}[{ \vert X_{2} \vert }^{\alpha_{2}}]{\bf E}[{ \vert X_3 \vert }^{\alpha_3}]
		({\rm Var}(Y_2))^{\frac{\alpha_2}{2}}({\rm Var}(Y_3))^{\frac{\alpha_3}{2}}
\,F\left(-\frac{\alpha_2}{2},-\frac{\alpha_3}{2};\frac{1}{2};\rho_{t}^{2}\right)\nonumber\\
		&\leq&{\bf E}[{ \vert X_{2} \vert }^{\alpha_{2}}]{\bf E}[{ \vert X_3 \vert }^{\alpha_3}]\,F\left(-\frac{\alpha_2}{2},-\frac{\alpha_3}{2};\frac{1}{2};\rho_{t}^{2}\right),
	\end{eqnarray}
where $ \rho_{t} \in [-1,1] $ is  the correlation coefficient
of $ Y_2 $ and $ Y_3 $
and $ F(\cdot) $ is a hypergeometric function defined by
	\begin{eqnarray*}\label{GHF}
		F(a,b;c;z):=\sum\limits_{n=0}^{+\infty}\frac{(a)_n(b)_n}{(c)_n}\cdot\frac{z^n}{n!},\  |z|\leq 1,
	\end{eqnarray*}
	with
	\begin{align*}
		(\alpha)_n:=\left\{
		\begin{array}{cl}
			\alpha(\alpha+1)\dots(\alpha+n-1), & n  \geq  1,\\
			1, & n=0,\alpha\neq 0.\\
		\end{array} \right.
	\end{align*}

The following property of the Gaussian hypergeometric function (see Andrews et al. \cite[(2.5.1)]{AAR99})
	\begin{align*}
		\frac{d}{dz}F(a,b;c;z)=\frac{ab}{c}F(a+1,b+1;c+1;z),
	\end{align*}
	and the Euler transformation (see  Andrews et al. \cite[Theorem 2.2.5]{AAR99} or Rainville \cite[Chapter 4, Theorem 21]{R60}) imply that
	\begin{align*}
		\frac{d}{dh}F\left(-\frac{\alpha_2}{2},-\frac{\alpha_3}{2};\frac{1}{2};h\right)
		&=\frac{\alpha_2\alpha_3}{2}F\left(1-\frac{\alpha_2}{2},1-\frac{\alpha_3}{2};\frac{3}{2};h\right)\\	
		&=\frac{\alpha_2\alpha_3}{2}(1-h)^{\frac{\alpha_2+\alpha_3}{2}-\frac{1}{2}}F\left(\frac{1}{2}+\frac{\alpha_2}{2},\frac{1}{2}+\frac{\alpha_3}{2};\frac{3}{2};h\right).
	\end{align*}
Then, when $ h \in [0,1] $,  $ \frac{d}{dh}F\left(-\frac{\alpha_2}{2},-\frac{\alpha_3}{2};\frac{1}{2};h\right)\ge 0. $  Therefore,  $ F\left(-\frac{\alpha_2}{2},-\frac{\alpha_3}{2};\frac{1}{2};h\right) $ reaches its maximum  value at $ h=1 $, and the maximum value is
	\begin{align*}
		\textbf{C}_{\alpha_2,\alpha_3}
		:=F\left(-\frac{\alpha_2}{2},-\frac{\alpha_3}{2};\frac{1}{2};1\right)
		=\frac{\Gamma(\frac{1}{2})\Gamma(\frac{1}{2}+\frac{\alpha_2+\alpha_3}{2})}
{\Gamma(\frac{1}{2}+\frac{\alpha_2}{2})\Gamma(\frac{1}{2}+\frac{\alpha_3}{2})},
	\end{align*}
which together with (\ref{proof-th-2-a}) and (\ref{proof-th-2-b}) implies that
	\begin{align*}
		{\bf E}[{ \vert X_1 \vert }^{-\alpha_1}{ \vert X_{2} \vert }^{\alpha_{2}}{ \vert X_3 \vert }^{\alpha_3}]
		&\leq\frac{{\bf E}[{ \vert X_{2} \vert }^{\alpha_{2}}]
{\bf E}[{ \vert X_3 \vert }^{\alpha_3}]}{\Gamma(1+\frac{\alpha_{1}}{2})}\int_{0}^{\infty}
\left(1+2t^{\frac{2}{\alpha_1}}\right)^{-\frac{1}{2}}
F\left(-\frac{\alpha_2}{2},-\frac{\alpha_3}{2};\frac{1}{2};\rho_{t}^{2}\right)dt\\
		&\leq\frac{{\bf E}[{ \vert X_{2} \vert }^{\alpha_{2}}]{\bf E}[{ \vert X_3 \vert }^{\alpha_3}]}{\Gamma(1+\frac{\alpha_{1}}{2})}\int_{0}^{\infty}
\left(1+2t^{\frac{2}{\alpha_1}}\right)^{-\frac{1}{2}}\textbf{C}_{\alpha_2,\alpha_3}dt\\
		&=\textbf{C}_{\alpha_2,\alpha_3}{\bf E}[{ \vert X_1 \vert }^{-\alpha_1}]{\bf E}
[{ \vert X_{2} \vert }^{\alpha_{2}}]{\bf E}[{ \vert X_3 \vert }^{\alpha_3}].
	\end{align*}
The proof is complete.

\subsection{Proof of Proposition  \ref{th-5}}

 The  inequality  on the right hand side of  \eqref{th-5-a} is a direct consequence of  Theorem \ref{th-1} when $n=3$.
 We now prove the inequality  on the left side of \eqref{th-5-a}.

For $t_1, t_2>0$, let $T:={\rm diag}(t_1^{\frac{2}{\alpha_1}},t_{2}^{\frac{2}{\alpha_{2}}},0)$. Then $\Sigma^{-1}+2T$ is positive
definite. Let $(Y_1, Y_2, Y_3)$
be a
%ZH centerd
centered
$\mathbb{R}^3$-valued Gaussian random vector with covariance matrix $( \Sigma^{-1}+2T)^{-1}$ and let
\begin{align*}
		\hat{f}(x_3):=\int_{\mathbb{R}^2}\frac{\exp(-\frac{1}{2}\left\langle \textbf{x},(\Sigma^{-1}+2T)\textbf{x}\right\rangle)}{(2\pi)^{\frac{n}{2}}
\sqrt{\det((\Sigma^{-1}+2T)^{-1})}}dx_1dx_2,
	\end{align*}
where  $\textbf{x}=(x_1, x_2, x_3)$.

Similar to the proof of Theorem \ref{th-1}, we  have that
	\begin{align*}
		{\bf E}[{ \vert X_1 \vert }^{-\alpha_1}{ \vert X_{2} \vert }^{-\alpha_{2}}{ \vert X_3 \vert }^{\alpha_3}]
		=\frac{{\bf E}[{ \vert X_3 \vert }^{\alpha_3}]}{\prod_{i=1}^{2}\Gamma(1+\frac{\alpha_{i}}{2})}\int_{0}^{\infty}
\int_{0}^{\infty}[\det(I_3+2T\Sigma)]^{-\frac{1}{2}}
({\rm Var}(Y_3))^{\frac{\alpha_3}{2}}dt_1dt_2,
	\end{align*}
and
	\begin{align*}
		{\rm Var}(Y_3)
		=\left(1+2{\bf t}^{'}[T_{1}^{-1}+2(\Sigma_{1}-{\bf tt}^{'})]^{-1}{\bf t}\right)^{-1},
	\end{align*}
	where $\Sigma_{1}:=\begin{pmatrix}
		1&a\\
		a&1
	\end{pmatrix}$, \   $T_{1}^{-1}:={\rm diag}(t_1^{-\frac{2}{\alpha_1}},t_{2}^{-\frac{2}{\alpha_{2}}}) $,\   $ {\bf t}:=\begin{pmatrix}
		b\\c
	\end{pmatrix} $.

Define
	\begin{align*}
		M:= [T_{1}^{-1}+2(\Sigma_{1}-{\bf tt}^{'})]
		=\begin{pmatrix}
			t_1^{-\frac{2}{\alpha_1}}+2(1-b^{2}) & 2(a-bc)\\
			2(a-bc) & t_{2}^{-\frac{2}{\alpha_{2}}}+2(1-c^{2})
		\end{pmatrix}.
	\end{align*}
Then	$ 	\det(M)=t_1^{-\frac{2}{\alpha_1}}t_{2}^{-\frac{2}{\alpha_{2}}}
+2(1-b^{2})t_{2}^{-\frac{2}{\alpha_{2}}}+2(1-c^{2})t_1^{-\frac{2}{\alpha_1}}
+4(1-a^{2}-b^{2}-c^{2}+2abc)$, and
	\begin{align*}
		M^{-1}=\frac{1}{\det(M)}\begin{pmatrix}
			t_{2}^{-\frac{2}{\alpha_{2}}}+2(1-c^{2}) &2(bc-a)\\
			2(bc-a) & t_1^{-\frac{2}{\alpha_1}}+2(1-b^{2})
		\end{pmatrix}.
	\end{align*}
Thus
	\begin{align*}
		{\bf t}^{'}M^{-1}{\bf t}&=\frac{b^2t_{2}^{-\frac{2}{\alpha_{2}}}+c^{2}t_1^{-\frac{2}{\alpha_1}}
+2b^2+2c^2-4abc}{t_1^{-\frac{2}{\alpha_1}}t_{2}^{-\frac{2}{\alpha_{2}}}
+2(1-b^{2})t_{2}^{-\frac{2}{\alpha_{2}}}+2(1-c^{2})t_1^{-\frac{2}{\alpha_1}}
+4(1-a^{2}-b^{2}-c^{2}+2abc)}\\
		&=\frac{b^2t_{1}^{\frac{2}{\alpha_{1}}}+c^{2}t_2^{\frac{2}{\alpha_2}}+
Kt_1^{\frac{2}{\alpha_1}}t_{2}^{\frac{2}{\alpha_{2}}}}{1+2(1-b^{2})
t_{1}^{\frac{2}{\alpha_{1}}}+2(1-c^{2})t_2^{\frac{2}{\alpha_2}}
+4\det(\Sigma)t_1^{\frac{2}{\alpha_1}}t_{2}^{\frac{2}{\alpha_{2}}}},
	\end{align*}
	where  $t_{1}, t_{2} \in (0,\infty)$ and $ K:=2b^2+2c^2-4abc$.
	
	First, we fix $ t_{1} $ below and calculate the maximum value of $ {\bf t}^{'}M^{-1}{\bf t} $ with respect to $ t_{2}$. We have
\begin{align*}
\frac{\partial}{\partial t_{2}}{\bf t}^{'}M^{-1}{\bf t}
		&=\frac{\frac{2c^2}{\alpha_2}t_2^{\frac{2}{\alpha_2}-1}+\left[
\frac{2K-4(1-c^2)b^2+4c^2(1-b^2)}{\alpha_2}\right]
t_1^{\frac{2}{\alpha_1}}t_2^{\frac{2}{\alpha_2}-1}+
\frac{4K(1-b^2)-8\det(\Sigma)b^2}{\alpha_2}t_1^{\frac{4}{\alpha_1}}
t_2^{\frac{2}{\alpha_2}-1}}
		{\left[1+2(1-b^{2})t_{1}^{\frac{2}{\alpha_{1}}}+2(1-c^{2})t_2^{\frac{2}{\alpha_2}}+
4\det(\Sigma)t_1^{\frac{2}{\alpha_1}}t_{2}^{\frac{2}{\alpha_{2}}}\right]^{2}}\\
		&:=\frac{\frac{2c^2}{\alpha_2}t_2^{\frac{2}{\alpha_2}-1}+\frac{I_1}
{\alpha_2}t_1^{\frac{2}{\alpha_1}}t_2^{\frac{2}{\alpha_2}-1}+
\frac{I_2}{\alpha_2}t_1^{\frac{4}{\alpha_1}}t_2^{\frac{2}{\alpha_2}-1}}
		{\left[1+2(1-b^{2})t_{1}^{\frac{2}{\alpha_{1}}}+2(1-c^{2})t_2^{\frac{2}{\alpha_2}}+
4\det(\Sigma)t_1^{\frac{2}{\alpha_1}}t_{2}^{\frac{2}{\alpha_{2}}}\right]^{2}},
	\end{align*}
	where $I_1:={2K-4(1-c^2)b^2+4c^2(1-b^2)}  $ and $ I_2:=4K(1-b^2)-8\det(\Sigma)b^2 $.
	By further computation, we  get that
	\begin{align*}
		I_1&={2K-4(1-c^2)b^2+4c^2(1-b^2)}\\
		&=2(2b^2+2c^2-4abc)-4b^2+4c^2\\
		&=8c^2-8abc,
	\end{align*}
and
	\begin{align*}
		I_2&=4K(1-b^2)-8\det(\Sigma)b^2\\
		&=4(2b^2+2c^2-4abc)(1-b^2)-8\det(\Sigma)b^2\\
		&=8[(1-a^2)-\det(\Sigma)](1-b^2)-\det(\Sigma)b^2\\
		&=8(1-a^2-b^2+a^2b^2)-8\det(\Sigma)\\
		&=8(ab-c)^2.
	\end{align*}
	Then, the numerator of $ \frac{\partial}{\partial t_{2}}{\bf t}^{'}M^{-1}{\bf t} $ is
\begin{align*}
		&\left[2c^2+(8c^2-8abc)t_1^{\frac{2}{\alpha_1}}+8(ab-c)^2t_1^{\frac{4}{\alpha_1}}\right]\frac{t_2^{\frac{2}{\alpha_2}-1}}{\alpha_2}\\
		&=2\left[\left(1+4t_1^{\frac{2}{\alpha_1}}+4t_1^{\frac{4}{\alpha_1}}\right)c^2+
\left(-4abt_1^{\frac{2}{\alpha_1}}-8abt_1^{\frac{4}{\alpha_1}}\right)c+
4a^2b^2t_1^{\frac{4}{\alpha_1}}\right]
\frac{t_2^{\frac{2}{\alpha_2}-1}}{\alpha_2}.
	\end{align*}
	Note that $ \left[\left(1+4t_1^{\frac{2}{\alpha_1}}+4t_1^{\frac{4}{\alpha_1}}\right)c^2+
\left(-4abt_1^{\frac{2}{\alpha_1}}
-8abt_1^{\frac{4}{\alpha_1}}\right)c+4a^2b^2t_1^{\frac{4}{\alpha_1}}\right]$  is a quadratic function of $ c $, and
	\begin{align*}
		\Delta&:=(-4abt_1^{\frac{2}{\alpha_1}}-8abt_1^{\frac{4}{\alpha_1}})^{2}-4(1+4t_1^{\frac{2}{\alpha_1}}+4t_1^{\frac{4}{\alpha_1}})4a^2b^2t_1^{\frac{4}{\alpha_1}}=0.
	\end{align*}
	Then we get that  $ \frac{\partial}{\partial t_{2}}{\bf t}^{'}M^{-1}{\bf t} \geq 0$. Therefore, for any fixed $ t_{1}$, the supermum of $ {\bf t}^{'}M^{-1}{\bf t} $ with respect to  $ t_{2} $  is
	\begin{align*}
		g(t_1)&:=\lim_{t_{2} \to \infty}{\bf t}^{'}M^{-1}{\bf t}\\
&=\lim_{t_{2} \to \infty}
		\frac{b^2t_{1}^{\frac{2}{\alpha_{1}}}+c^{2}t_2^{\frac{2}{\alpha_2}}+Kt_1^{\frac{2}{\alpha_1}}t_{2}^{\frac{2}{\alpha_{2}}}}{1+2(1-b^{2})t_{1}^{\frac{2}{\alpha_{1}}}+2(1-c^{2})t_2^{\frac{2}{\alpha_2}}+4\det(\Sigma)t_1^{\frac{2}{\alpha_1}}t_{2}^{\frac{2}{\alpha_{2}}}} \\
		&=\frac{c^2+Kt_1^{\frac{2}{\alpha_1}}}{2(1-c^2)+4\det(\Sigma)t_1^{\frac{2}{\alpha_1}}}.
	\end{align*}

Now, we take the derivative of $ g(t_1) $ with respect to $ t_1$ and get that
	\begin{align*}
		\frac{d}{dt_1}g(t_1)&=\frac{\frac{4K(1-c^2)}{\alpha_1}t_1^{\frac{2}{\alpha_1}-1}-
\frac{8\det(\Sigma)c^2}{\alpha_1}t_1^{\frac{2}{\alpha_1}-1}}
{\left[2(1-c^2)+4\det(\Sigma)t_1^{\frac{2}{\alpha_1}}\right]^2}.
	\end{align*}
	In the above expression,  we know that   $ \frac{d}{dt_1}g(t_1)\ge0  $  if and only if $ 4K(1-c^2)-8\det(\Sigma)c^2 \ge 0.$
However,
	\begin{align*}
		4K(1-c^2)-8\det(\Sigma)c^2&=4[(2b^2+2c^2-4abc)(1-c^2)-2\det(\Sigma)c^2]\\
		&=8[[(1-a^2)-\det(\Sigma)](1-c^2)-\det(\Sigma)c^2]\\
		&=8(1-a^2-c^2+a^2c^2-\det(\Sigma))\\
		&=8(ac-b)^2\geq0.
	\end{align*}
	Therefore, $ \frac{d}{dt_1}g(t_1)\geq0 $. Then, the supermum of $ g(t_1) $ is
	\begin{align*}
		\lim_{t_{1} \to \infty}g(t_1)
		&=\lim_{t_{1} \to \infty}
		\frac{c^2+Kt_1^{\frac{2}{\alpha_1}}}{2(1-c^2)+4\det(\Sigma)t_1^{\frac{2}{\alpha_1}}}\\
		&=\frac{K}{4\det(\Sigma)}=\frac{(1-a^2)-\det(\Sigma)}{2\det(\Sigma)}\\
&=\frac{(1-a^2)}{2\det(\Sigma)}-\frac{1}{2}.
	\end{align*}
	Then
	\begin{align*}
		{\rm Var}(Y_3)
		=\left(1+2{\bf t}^{'}[T_{1}^{-1}+2(\Sigma_{1}-{\bf tt}^{'})]^{-1}{\bf t}\right)^{-1}\geq\frac{\det(\Sigma)}{(1-a^2)}.
	\end{align*}
Hence, by virtue of (\ref{a3}), we obtain
	\begin{align*}
		{\bf E}\left[{ \vert X_1 \vert }^{-\alpha_1}{ \vert X_{2} \vert }^{-\alpha_{2}}{ \vert X_3 \vert }^{\alpha_3}\right]
		\geq\frac{\det(\Sigma)}{1-a^{2}}{\bf E}\left[{ \vert X_1 \vert }^{-\alpha_1}\right]{\bf E}\left[{ \vert X_{2} \vert }^{-\alpha_{2}}\right]{\bf E}\left[{ \vert X_3 \vert }^{\alpha_3}\right].
	\end{align*}
The proof is completed.

\bigskip

{ \noindent {\bf\large Acknowledgments}  This work was supported
by the Scientific  Foundation of  Nanjing University of Posts and Telecommunications (NY221026),  the National Natural Science Foundation of China
%RS (12171335)
(12171335, 12071011, 11931004)
and the Science Development Project of Sichuan University (2020SCUNL201).

\end{document}